\documentclass[centertags,leqno]{article}

\usepackage{amsmath}
\usepackage{amssymb}
\usepackage{latexsym}
\usepackage{amsxtra}
\usepackage{amscd}
\usepackage{theorem}
\usepackage{xypic}

\theoremstyle{plain}

{\theorembodyfont{\slshape}

        \newtheorem{thm}{Theorem}[section]
        
        \newtheorem{lem}[thm]{Lemma}
        \newtheorem{prop}[thm]{Proposition}

}

{\theorembodyfont{\rmfamily}

        \newtheorem{defn}[thm]{Definition}
        
        \newtheorem{rem}[thm]{Remark}

        \newtheorem{exa}[thm]{Example}

}

\renewcommand{\em}{\sl}

\newcommand{\proof}{{\bf Proof:\ }}

\newcommand{\Endproof}{\hspace*{\fill} $\Box$ \vspace{1ex} \noindent }

\makeatletter
\renewcommand{\subsection}{\@startsection{subsection}{2}%
        {\z@}{-3.25ex plus -1ex minus-.2ex}{-1em}{\bf}}
\makeatother

\newcommand{\PP}{\mathbb{P}}
\newcommand{\ZZ}{\mathbb{Z}}

\newcommand{\EE}{\mathbb E}

\newcommand{\FF}{\mathbb{F}}
\newcommand{\A}{\mathbb{A}}

\newcommand{\BB}{\mathbb{B}}

\newcommand{\C}{\mathcal{C}}

\newcommand{\T}{\mathcal{T}}
\newcommand{\E}{\mathcal{E}}
\newcommand{\M}{\mathcal{M}}
\newcommand{\N}{\mathcal{N}}
\newcommand{\LL}{\mathcal{L}}
\newcommand{\OO}{\mathcal{O}}

\newcommand{\Zb}{\bar{Z}}
\newcommand{\G}{\mathcal{G}} 
\newcommand{\lt}{{\rm\scriptstyle lt}}
\newcommand{\I}{\mathcal{I}}
\newcommand{\olog}{\Omega^{\rm\scriptscriptstyle log}}

\newcommand{\Hom}{\mathop{\rm Hom}\nolimits}
\newcommand{\End}{\mathop{\rm End}\nolimits}

\newcommand{\ord}{{\rm ord}}

\newcommand{\mup}{\boldsymbol{\mu}_p}

\newcommand{\Res}{\mathop{\rm Res}\nolimits}

\newcommand{\To}{\longrightarrow}

\newcommand{\Spec}{\mathop{\rm Spec}\nolimits}

\newcommand{\Cf}{{\mathfrak C}}
\newcommand{\Def}{\mathop{\rm Def}}
\newcommand{\bsigma}{{\boldsymbol\sigma}}
\newcommand{\PGL}{\mathop{\rm PGL}\nolimits}

\title{The accessory parameter problem in positive characteristic}
\author{Irene I.\ Bouw}
\date{}
\begin{document}
\maketitle
\begin{abstract} We study the existence of Fuchsian differential equations in 
positive characteristic with nilpotent $p$-curvature, and given local
invariants. In the case of differential equations with logarithmic
local mononodromy, we determine the minimal possible degree of a
polynomial solution.\\[2ex]
  2000 Mathematical Subject Classification: Primary  14D10, 12H20
\end{abstract}

This paper deals with second order differential equations with regular
singularities in characteristic $p>0$. Our main interest is to
characterize those differential equations with nilpotent (resp.\
nilpotent but nonzero) $p$-curvature. This problem is known as Dwork's
{\sl accessory parameter problem}.  Differential equations with
nilpotent $p$-curvature arise naturally in algebraic geometry. For
example, Katz (\cite{Katz}) showed that differential equations
``coming from geometry'', such as Picard--Fuchs differential
equations, have nilpotent $p$-curvature.  The nilpotence (resp.\
nonvanishing) of the $p$-curvature may be characterized in terms of
the existence of polynomial solutions.  The study of polynomial
solution of differential equations in positive characteristic goes
back to Dwork (\cite{Dwork}, \cite{Dwork2}) and Honda (\cite{Honda}).

More recently, differential equations with nilpotent but nonzero $p$-curvature
came up in Mochizuki's work on $p$-adic uniformization
(\cite{Mochizuki1}, \cite{Mochizuki2}). Mochizuki develops the theory
of indigenous bundles.  On a curve of genus zero, these
may be identified with differential equations with nilpotent but
nonzero $p$-curvature (\cite[\S 5]{indi}). To prove concrete existence
results, the description of indigenous bundles as differential
equations turns out to be more convenient.

Differential equation with nilpotent but nonzero $p$-curvature also
arise in the theory of reduction to characteristic $p$ of Galois
covers of curves.  Solutions of differential equations arise in this
context in the form of {\sl deformation data}. In \cite{indi} one
finds a correspondence between indigenous bundles and deformation
data. (See also \S \ref{corrsec}.)

This paper combines results from the work of Mochizuki on indigenous
bundles with results on deformation data and techniques from the work
of Honda and Dwork.  It turns out that combing
these techniques is very fruitful and allows to answer 
questions which are interesting from all three points of view.

\bigskip\noindent
We now give a more detailed description of the content of the paper.

 We fix the number, $r$, of singularities of the differential
equation, together with a set ${\boldsymbol \alpha}=(\alpha_1, \ldots,
\alpha_r)$ of local invariants, the {\sl local exponents}. We study
the stack $\N_{0,r}({\boldsymbol \alpha}; n)$ of differential
operators, $L$, with local exponents ${\boldsymbol \alpha}$, nilpotent
but nonzero $p$-curvature, and  strength $n$. The {\sl strength} is a
natural invariant of the differential operator introduced by Mochizuki
(\cite{Mochizuki1}) which is defined as the degree of the zero divisor of the
$p$-curvature (\S \ref{accsec}). Essentially, it corresponds to the
minimal degree of a polynomial solution of $L$. We refer to \S
\ref{accsec} for a precise statement.

The main question we are interested in is to determine (${\boldsymbol
\alpha}, n$) such that $\N_{0,r}({\boldsymbol \alpha}; n)$ is nonempty, 
and to determine the dimension of the irreducible components of those
$\N_{0,r}({\boldsymbol \alpha}; n)$. This variant of Dwork's accessory
parameter problem  we call the {\sl strong accessory parameter
problem}.

Our strongest results are in the case of {\sl logarithmic local
monodromy} (\S \ref{logmonsec}) which is the case on which the work of
Mochizuki focusses.  Mochizuki shows that if $\N_{0,r}({\boldsymbol
\alpha}; n)$ is nonempty, then every irreducible component has the
maximal possible dimension $r-3$.  We determine all those $n$ for
which $\N_{0,r}({\boldsymbol \alpha}; n)$ is nonempty. Equivalently,
we determine which degrees occur as the minimal degree of a polynomial
solution of a differential operator $L$.

More precisely, we show the following.

\bigskip\noindent{\bf Theorem \ref{nonemptythm}} Suppose that
  $p>r-2$. Then the stack $\N_{0,r}({\boldsymbol \alpha}=(0,\ldots,
  0); n)$ is nonempty if and only if the strength, $n$, is congruent
  to $0\bmod{2p}$ and satisfies $0\leq n\leq (r-2)(p-1)/2$.

\bigskip\noindent The necessity of these conditions already follows
from a result of Dwork (\cite[Lemma 10.1]{Dwork2}. To prove the
theorem, we first apply a deformation technique due to Mochizuki to
reduce ourselves to the case that $r=n/p+3$. In this case, the minimal
degree, $d$, of a polynomial solution is less than $p$. We then
explicitly construct differential equations with a solution of this  degree
$d$ (Proposition \ref{recursionprop}). Our method here is inspired by
the work of Honda, Dwork and Beukers.

For general choice of the local exponents ${\boldsymbol \alpha}$ our
results are less strong.  We give a new necessary condition on $n$ for
$\N_{0, r}({\boldsymbol \alpha}; n)$ to be nonempty (Lemma
\ref{nlem}). Our main result in this case is a result on the dimension
of the irreducible components of $\N_{0, r}({\boldsymbol \alpha}; n)$
(Proposition \ref{dimmaxprop}). This is a (weaker) analog of the
result which Mochizuki proved in the case of logarithmic local monodromy. 

The
key ingredient in the proof of Proposition \ref{dimmaxprop} is the
deformation theory of deformation data (\S \ref{mupsec}), following
Wewers (\cite{cotang}). We use this as a replacement for Mochizuki's
results on deformation of indigenous bundles which are not available
here.  We also give some concrete examples which illustrates what to
expect in the general case (\S \ref{exasec}).

\section{Fuchsian differential equations in positive characteristic}
\label{Fuchssec}
The main goal of this section is to recall and reformulate some
classical results on differential equations in positive
characteristic. In \S \ref{polysolsec} we recall some results of
Dwork and Honda on polynomial solutions of differential equations. In
\S \ref{pcurvsec} we recall the definition of the $p$-curvature and
state some basic properties.

\subsection{Algebraic solutions of Fuchsian differential equations}
\label{polysolsec}
Let $X=\PP^1_k$, and choose a parameter $x$ on $X$. Let $r\geq 3$, and
suppose given pairwise distinct points $x_1, \ldots, x_r\in X$. We
assume that $x_r=\infty$.

In this paper we consider {\sl Fuchsian differential operators}
\begin{equation}\label{Leq}
L=(\partial/\partial x)^2+p_1(\partial/\partial x)+p_2
\end{equation}
on $X$ with singularities in $x_1, \ldots, x_r=\infty$. Recall that
this means that $p_1$ and $p_2$ are rational functions on $X$ which
are regular outside the $x_i$ such that $\ord_{x_j} p_i\geq i$ for all
$j$. We call $L(u)=0$ the corresponding differential equation.  The
differential operator  $L_W=(\partial/\partial x)+p_1$ is called the
{\sl Wronskian equation} associated to $L$ (see \cite[\S 1]{Honda} for
the relation of $L_W$ with the Wronskian.) All differential operators
in this paper  have order $2$ and are supposed to be Fuchsian.

\begin{defn}\label{pcurvdef}
Let $L$ be as in (\ref{Leq}).
\begin{itemize}
\item[(a)] We say that $L$ has {\sl nilpotent $p$-curvature} if both
$L$ and $L_W$ have a polynomial solution.
\item[(b)] Let $L$ be a differential operator with nilpotent
$p$-curvature. We say that $L$ has {\sl nonzero $p$-curvature} if
the space of polynomial solutions of $L$ is $1$-dimensional over $k[x]^p$.
\end{itemize}
\end{defn}

We refer to \cite{Honda} for a discussion of these notions. Nilpotent
$p$-curvature is called ``sufficiently many solutions in a weak
sense'' by Honda (\cite{Honda}). To $L$ we may associate a flat vector
bundle $(\E, \nabla)$ of rank $2$ on $X$. To this flat vector bundle
is naturally associated the $p$-curvature (see for example
\cite{Katz}). In \cite[Appendix]{Honda} it is shown that the
$p$-curvature of $\E$ is nilpotent (resp.\ nonzero) if and only if $L$
satisfies the conditions of Definition \ref{pcurvdef}. For a detailed
discussion of the correspondence of $\E$ and $L$ we refer to \cite[\S
5]{indi}. In \S \ref{pcurvsec} we give a short introduction. For a
definition of the $p$-curvature in terms of $L$ we refer to
\cite{Dwork}.

Since $L$ has regular singular points in $x_1, \ldots, x_r=\infty$, we
may write
\[
p_1= \frac{P_1}{P_0}, \qquad p_2=\frac{P_2}{P_0}+\frac{P_3}{P_0^2},
\]
where $P_0=\prod_{i=1}^{r-1}(x-x_i).$ For $i=1, \ldots r-1$, we define
the {\sl local exponents} $\{\alpha_i, \alpha'_i\}\in \FF_p^\times$ of
$L$ at $x_i$ by
\begin{equation}\label{locexpeq}
\alpha_i+\alpha'_i=1-\frac{P_1(x_i)}{P_0'(x_i)}, \qquad
\alpha_i\alpha'_i=\frac{P_3(x_i)}{P_0'(x_i)^2}.  
\end{equation}
Similarly, for $x=x_r=\infty$ we define the {\sl local exponents}
$\{\alpha_r, \alpha'_r\}\in \FF_p^\times$ of $L$ at $x_r$ by
\begin{equation}\label{locexpeq2}
\alpha_r+\alpha'_r=\left(-1+\frac{P_1(x)}{x^{r-2}}\right)({x=\infty}), \qquad
\alpha_r\alpha'_r=\left(\frac{P_2(x)}{x^{r-3}}\right)({x=\infty}).
\end{equation}
The local exponents are the eigenvalues of the local monodromy matrix of $L$ at
$x=x_i$.  They satisfy the Riemann relation
\begin{equation}\label{Rreq}
\sum_{i=1}^r (\alpha_i+\alpha_i')=r-2.
\end{equation}
This follows by immediate verification (cf.\ \cite{Yoshida}).

\begin{lem}\label{locexplem}
\begin{itemize}
\item[(a)]
Suppose that $L$ has nilpotent $p$-curvature. Then the local exponents
$\alpha_i, \alpha'_i$ of $L$ at $x_i$ are in $\FF_p$.
\item[(b)] Let $L$ be a differential operator as in (\ref{Leq}) whose
local exponents $\{\alpha_i, \alpha_i'\}$ are all elements of
$\FF_p$. Suppose that $L$ has a solution $u\in k[x]$. Then $L$ has
nilpotent $p$-curvature.
\end{itemize}
\end{lem}

\proof Part (a) is proved in \cite[Prop. 2.1]{Honda}.  If the local
exponents of $L$ are in $\FF_p$ then (\ref{denormaleq}) implies that
$p_1=Q'/Q$ for some polynomial $Q\in k[x]$. Part (b) immediately
follows from this observation (cf.\ \cite[Cor.\ 1 to Prop.\
2.3]{Honda}.)  \Endproof

\begin{lem}\label{deglem}
Let $u=u(x)$ be a polynomial solution of $L(u)=0$. Then
\begin{itemize}
\item[(a)] $\deg(u)\equiv -\alpha_r\pmod{p}$ or $\deg(u)\equiv
-\alpha_r'\pmod{p}$,
\item[(b)] $\ord_{x_i}(u)\equiv \alpha_i$ or $\ord_{x_i}(u)\equiv
\alpha_i'\pmod{p}$ for $i\neq r$.
\end{itemize}
\end{lem}

\proof Suppose that $x$ is a local parameter of $X=\PP^1_k$ at the
singular point $x_i$ for $1\leq i\leq r$. Rewriting the differential
equation in terms of $x$ immediately yields the statement of the
lemma; (a) corresponds to $i=r$ and (b) corresponds to $i\neq r$.
\Endproof

The following proposition is proved by Honda. It is a stronger version
of Lemma \ref{deglem}.(b) in the case that $L$ has zero $p$-curvature.

\begin{prop}[Honda]\label{pcurv0prop}
Let $L$ be a differential operator whose $p$-curvature is zero. Let
$u_1, u_2\in k[x]$ be two solutions of $L$ which are independent over
$k[x]^p$. Suppose, moreover, that $\delta:=\deg(u_1)+\deg(u_2)$ is
minimal.
\begin{itemize}
\item[(a)]
Then $\deg(u_1)$ and $\deg(u_2)$ are noncongruent to each
other (modulo $p$). In particular, $\{\deg(u_1), \deg(u_2)\}\equiv
\{-\gamma_1, -\gamma_2\} \pmod{p}$ and $\gamma_1\not\equiv
\gamma_2\pmod{p}$.
\item[(b)] Suppose that $p\neq 2$ and $p\geq r-2$. Then 
\[ 
r\leq \delta\leq (r-1)p-(2r-3).
\]
\end{itemize}
\end{prop}

\proof
 Part (a) is \cite[Prop.\
5.1]{Honda}. Part (b) is \cite[Theorem 7.c]{Honda}
\Endproof

We conclude from Proposition \ref{pcurv0prop} that if $L$ has
nilpotent $p$-curvature then there exists a unique monic polynomial
solution of minimal degree.

Let $L$ be a differential operator as above with nilpotent
$p$-curvature. Let $u$ be a polynomial solution of minimal
degree. Then $\tilde{\alpha}_i:=\ord_{x_i} u <p$ for $i\neq
r$. We may suppose that $\tilde{\alpha}_i\equiv \alpha'_i
\pmod{p}$ for $i\neq r$ (Lemma \ref{deglem}.(b)). For $i=r$, we
suppose that $\deg(u)\equiv -\alpha'_r\pmod{p}$.  It follows that we
may write
\[
u=\prod_{i=1}^{r-1}(x-x_i)^{\tilde{\alpha}_i} v, \qquad v\in k[x].
\]

 For $i=1, \ldots, r$, we define integers $0\leq t_i\leq p-1$ by
\begin{equation}\label{tieq}
t_i\equiv \alpha_i'-\alpha_i\pmod{p}.
\end{equation}

The following proposition is an analog of Proposition \ref{pcurv0prop}.(b) in 
the case that the $p$-curvature is nonzero.

\begin{prop}[Dwork] \label{Dworkdegprop} Let $L$ be a differential operator
 with 
nilpotent but nonzero $p$-curvature. Let $u$ be a polynomial solution
of minimal degree.  There are exists a nonnegative integer $t$ such
that
\[
2\deg(v)+pt=(p-1)(r-2)-(\sum_{i=1}^{r}t_i).
\]
\end{prop}

\proof
This is \cite[Lemma 10.1]{Dwork2}
\Endproof

The following lemma is an easy consequence of Lemma \ref{deglem}. If
$L$ is a differential operator with nilpotent $p$-curvature, the the
lemma gives a sufficient criterion for the $p$-curvature of $L$ to be
nonzero.

\begin{lem}\label{pcurvnonzerolem}
Let $L$ be a differential operator with nilpotent
$p$-curvature. Suppose that $\alpha_i=\alpha_i'$, for some $i$. Then
the $p$-curvature is nonzero.
\end{lem}

\proof Let $L$ be as in the statement of the lemma, and suppose that
$\alpha+i=\alpha_i'$.  In the case that $i=r$, the lemma follows from
Lemma \ref{deglem}.(a). The general case is easily reduced to this
case, by changing the coordinate $x$.  
\Endproof

Beukers (\cite{Beukers}) proves the same result in the case of $r=4$
singularities such that $\alpha_i=\alpha_i'=0$ for all $i$.

\subsection{Normalized differential operators}\label{normalsec}
\begin{defn}\label{equivdef}
Two differential equations $L_1$ and $L_2$ are called {\sl equivalent}
if they have the same set of singularities and there exists a rational
function $v$ on $X$ such that $L_1(u)=0$ if and only if $L_2(vu)=0$
for all $u\in k(\!(x)\!)$.
\end{defn}

Let $L_1$ and $L_2$ be equivalent differential operators, and let
$\{\alpha_i, \alpha'_i\}$ be the local exponents of $L_1$ at
$x_i\neq \infty$. There exist  $\mu_i\in \FF_p$ such that the local
exponents of $L_2$ at $x_i$ are $\{\alpha_i-\mu_i,
\alpha'_i-\mu_i\}$. Here $\mu_i=\ord_{x_i} v$, for $v$ as in
Definition \ref{equivdef}. 

\begin{defn}\label{udef}\label{normaldef}
A differential equation $L$ with nilpotent $p$-curvature is called
{\sl normalized} if its (unique) monic solution of minimal degree does
not have zeros in the singular points $x_i$.  We denote the monic
solution of minimal degree by $u$, and put $d=\deg(u)$.
\end{defn}

The local exponents of a normalized differential
operator are uniquely determined by ${\boldsymbol \alpha}:=\{\alpha_1,
\ldots, \alpha_r\}$, due to the Riemann relation (\ref{Rreq}). Namely,
the local exponents are $(0, \alpha_i)$ at $x=x_i$ for $i\neq r$ and
$(-d, -d+\alpha_r)$ for $x=x_r=\infty$.

 Normalized differential operators satisfy
\begin{equation}\label{denormaleq}
p_1=\sum_{i=1}^{r-1} \frac{1-\alpha_i}{x-x_i}, \qquad
p_2=\frac{d(d-\alpha_r)x^{r-3}+\beta_{r-4}x^{r-4}+\cdots
+\beta_0}{\prod_{i=1}^{r-1}(x-x_i)}.  
\end{equation}

Note that the properties `nilpotent $p$-curvature' and `nonzero
$p$-curvature' only depends of the equivalence class of the
differential operator. Moreover, every equivalence class of
differential operators with nilpotent $p$-curvature contains a
normalized one.

\subsection{The $p$-curvature}\label{pcurvsec}
Let $L=(\partial/\partial x)^2+p_1(\partial/\partial x)+p_2$ be a
Fuchsian differential operator defined over $k$ with nilpotent
$p$-curvature. Let $u\in k[x]$ be a polynomial solution of minimal
degree. Lemma \ref{locexplem}.(a) implies that there exists a
polynomial $Q\in k[x]$ such that $p_1=Q'/Q$. An explicit formula for
$Q$ is easily deduced from (\ref{denormaleq}).

 Let $(L, u)$ be as above. Let $D:=\partial/\partial x$. We define a
flat vector bundle $\E$ on $\PP^1$, as in the proof of
\cite[Proposition 5.3]{indi}. Namely, on $\A^1$ we let $\E$ be the
trivial bundle with basis $e_1, e_2$ and connection defined by
\[
\nabla(D)
\begin{pmatrix}e_1\\e_2\end{pmatrix}=\begin{pmatrix} 0&-p_2\\
1&-p_1\end{pmatrix} \begin{pmatrix}e_1\\e_2\end{pmatrix}.
\]
It is shown in loc.cit.\ that $\E$ extends to a flat vector bundle
with regular singularities on $\PP^1$. 

Let $\T=(\olog_X/k)^{\otimes -p}$.
The $p$-curvature is an $\OO_X$-linear map
\[
\Psi_\E:\T\to \End_{\OO_X}(\E).
\]
Since $D^p=0$, the $p$-curvature of $\E$ is defined by
$\Psi_\E(D^{\otimes p}):=\nabla(D)^p$ (cf.\ \cite[\S 3.1]{indi}).  It
is well known that the $p$-curvature of $L$ is nilpotent (resp.\ zero)
if and only if the matrix of $\Psi_\E(D^{\otimes p})$ is nilpotent
(resp.\ zero) (\cite[Appendix]{Honda}).

\begin{prop}\label{pcurvprop}
The differential operator $L$ has nonvanishing $p$-curvature if and only if
\begin{equation}\label{pcurveq}
\left(\frac{\partial}{\partial x}\right)^{p-1}\frac{1}{Qu^2}\neq 0.
\end{equation}
\end{prop}

\proof
 Suppose that (\ref{pcurveq}) holds.  The computation in
the proof of \cite[Proposition 5.3]{indi} implies that
\[
\Psi_\E(D^{\otimes p})\frac{e_1}{u}\neq 0.
\]
We conclude that the $p$-curvature of $L$ is nonzero.

To prove the other implication, we suppose that $\Psi_\E(D^{\otimes
p})\neq 0$. Let $\M\subset \E$ be the kernel of $\Psi_\E$ (\cite[\S
3.1]{indi}). Our assumption implies that this is a flat subbundle of
$\E$ of rank $1$.  We may choose a horizontal section $\eta$ of $\M$
which generates $\M$ on a dense open subset of $\PP^1_k$.  Reversing
the computation from \cite[Proposition 5.3]{indi} implies that
(\ref{pcurveq}) holds.  \Endproof

\section{Deformation data}\label{ddgensec}
In this section, we prove a correspondence between deformation data
and differential equations with nilpotent but nonzero
$p$-curvature. This is a reformulation and adaption to the present
case of the result of \cite{indi}. For the convenience of the reader, we
start by recalling the definition and basic properties of deformation
data. We refer to \cite{special, bad} for more details.  A short
introduction explaining how deformation data arise naturally in the
theory of stable reduction can be found in \cite[\S 4.2]{indi}. 

\subsection{Definitions}\label{dddefsec}
Let $k$ be an algebraically closed field of characteristic $p>2$.

\begin{defn}\label{dddef}
  A {\sl deformation datum} of type $(H,\chi)$ is a pair $(g,
  \omega)$, where $g:Z\to X=\PP^1_k$ is a finite Galois cover of
  smooth projective curves and $\omega$ is a meromorphic differential
  form on $Z$ such that the following conditions hold.
\begin{itemize}
\item[(a)] 
Let $H$ be the Galois group of $Z\to X$.
Then
\[\beta^\ast \omega=\chi(\beta)\cdot
  \omega, \qquad \mbox{for all } \beta\in H.
\]
 Here $\chi:H\to
  \FF_p^\times$ in an injective character.
\item[(b)] The differential form $\omega$ is logarithmic, i.e.\ of the
  form $\omega={\rm d} f/f$, for some meromorphic function $f$ on
  $Z$.
\end{itemize}
\end{defn}

Part (b) of Definition \ref{dddef} may also be reformulated as:
$\omega$ is fixed by the Cartier operator $\C$.  Let $(g, \omega)$ be
a deformation datum.  For each closed point $x\in X$ we define the
following invariants.
\begin{equation}\label{sigmaeq}
m_x:=|H_z|, \qquad h_x:=\ord_z(\omega)+1, \qquad \sigma_x:=h_x/m_x.
\end{equation}
Here $z\in Z$ is some point above $x$ and $H_z\subset H$ is the
stabilizer of $z$. The invariant $\sigma_x$ is called the {\sl
ramification invariant} of the deformation datum at $x$.
 The
 following lemma gives some necessary conditions on the ramification
 invariant $\sigma_x$.

\begin{lem} \label{defodatlem}
  Let $(Z,\omega)$ be a deformation datum. 
  \begin{itemize}
  \item[(a)] For all $x\in X$, the ramification index $m_x$ divides
    $p-1$. Moreover, $h_x$ and $m_x$ are relatively prime.
  \item[(b)]
    If $h_x\not=0$ then $\gcd(p,h_x)=1$.
  \item[(c)]
    For all but finitely many points $x\in X$ we have $\sigma_x=1$.
  \item[(d)]
    We have
    \[
          \sum_{x\in X} (\sigma_x-1) = 2g-2.
    \]
  \end{itemize}
\end{lem}

\proof
This lemma is proved in \cite[Lemma 4.3]{indi}.
\Endproof

\begin{defn}\label{ddssdef}
Let $(Z,\omega)$ be a deformation datum on $X$. 
\begin{itemize}
\item[(a)] A point $x\in X$ is said to be a {\sl critical point} of
the deformation datum if $\sigma_x\neq 1$.
\item[(b)] A critical point $x\in X$
is   {\em supersingular }  if
$\sigma_x=(p+1)/(p-1)$. 
\item[(c)] A critical point $x\in X$ is  {\em singular}
if it is {\em not} a supersingular point and
$\sigma_x\not\equiv 1\pmod{p}$. 
\item[(d)] A critical point $x\in X$ such that $\sigma_x\equiv
1\pmod{p}$ is called a {\sl spike}.
\end{itemize}
\end{defn}

We refer to \cite{indi} for a motivation of the terminology.  Lemma
 \ref{defodatlem} (c) implies that a deformation datum has  finitely many
 critical points.

\begin{defn} \label{signaturedef}
  A {\em signature} is given by a finite set $M$ and a map
  \[
      \bsigma: M \to \frac{1}{p-1}\cdot\ZZ,\quad x \mapsto \sigma_x
  \]
  such that $\sigma_x\geq 0$ and $\sigma_x\neq 1,(p+1)/(p-1)$ for
  all $x\in M$ and such that the number
  \[
        d := \frac{p-1}{2}\Big(\,2g-2
        -\sum_{x\in M} (\sigma_x-1) \,\Big)
  \]
  is a nonnegative integer. The {\em singularities} of $\bsigma$ are the
  elements $x\in M$ with $\sigma_x\not\equiv 1\pmod{p}$. 

  Given a deformation datum $(Z,\omega)$ on $X$, the invariants
  $\sigma_x$ defined in \eqref{sigmaeq} give rise to a signature
  $\bsigma$ (where $M$ is the set of points $x\in X$ with
  $\sigma_x\neq 1,(p+1)/(p-1)$). It follows from Lemma
  \ref{defodatlem} (d) that the number $d$ defined above is the number of
  supersingular points.
\end{defn}

We denote by $\BB$ the set of critical points of the deformation
datum.  Let $r$ be the number of singularities of the deformation
datum.  We say that the deformation datum $(Z,\omega)$ is {\em
trivial} if $2g-2+r=0$.  In the rest of this paper, we exclude trivial
deformation data.  Without loss of generality, we may therefore
suppose that $x_r=\infty$ is a singularity. Let $\BB'=\{b\in \BB\,|\,
x_b\neq \infty\}$.

We let $s-r$ be the number of spikes.  We always
enumerate the singularities (resp.\ spikes) of a deformation datum
$(Z,\omega)$ as $x_1,\ldots,x_r$ (resp.\ $x_{r+1}, \ldots, x_s$) and
write $\sigma_i:=\sigma_{x_i}$ for $i=1,\ldots,s$.

  Define integers $0\leq a_i< p-1$ and $ \nu_i\geq 0$ by
\begin{equation}\label{aidefeq}
\sigma_i=\frac{a_i}{p-1}+\nu_i.
\end{equation}
Note that $x$ is supersingular if and only if $(a_x, \nu_x)=(2, 1)$.

Let $(Z, \omega)$ be a deformation datum of signature ${\boldsymbol
\sigma}=(\sigma_i)$.  Let $u\in k[x]$ be the monic polynomial whose
zeros are exactly the supersingular points of the deformation datum
(with multiplicity one).

For future reference we note that $Z$ is a connected component of the
smooth projective curve defined by the Kummer equation
\begin{equation}\label{Kummereq}
z^{p-1}=\prod_i (x-x_i)^{a_i}u^2,
\end{equation}
cf.\ (\ref{aidefeq}).  Let $S\subset Z$ be the inverse image of the
set of critical points $x_i\in X$ for which $\sigma_i=0$. The
definition of $\sigma_i$ implies that $S$ is the set of poles of
$\omega$. Therefore $\omega$ is a section of the sheaf
$\olog:=\Omega_{Z/k}(S)$ of differential $1$-forms with simple poles
in $S$. Definition \ref{dddef}.(a) implies therefore that $\omega\in
H^0(Z, \olog)_\chi$.

Once $Z$ and $\chi$ are given, we may characterize the deformation
data $(Z, \omega)$ as those section of $H^0(Z, \olog)_\chi$ which are
fixed by the Cartier operator $\C$. The following lemma is stated for
completeness.

\begin{lem}\label{logpolebasislem}
\begin{itemize}
\item[(a)] The dimension of $H^0(Z, \olog)_\chi$ is $|\BB|-1-
(\sum_i a_i)/(p-1)$.
\item[(b)]
The differentials
\[
\omega_j=\frac{x^{j-1}z{\rm d}x}{\prod_{i\in\BB'}(x-\tau_i)}, \qquad
j=1, \ldots ,r-1-(\sum_i a_i)/(p-1),
\]
form a basis of $H^0(Z, \olog)_{\chi}$.
\end{itemize}
\end{lem}

\proof This is proved like \cite[Lemma 4.3]{cyclic}.
\Endproof

\subsection{A correspondence}\label{corrsec}
In this section we prove a correspondence between deformation data and
solutions of Fuchsian differential equation in positive
characteristic.

Let $k$ be an algebraically closed field of characteristic $p$. As
before, we choose a parameter $x$ on $X=\PP^1_k$. We write
$D=\partial/\partial x$ and $f'$ for $D(f)$.

 Suppose we are given a deformation datum $(Z, \omega)$ of type
$(H,\chi)$ of signature ${\boldsymbol \sigma}$. As before we denote by
$x_1, \ldots, x_r=\infty$ the singular critical points of the
deformation datum and by $x_{r+1}, \ldots, x_s$ the spikes. After
replacing $k$ by a larger algebraically closed field, if necessary,
we may suppose that all $x_i$ are $k$-rational. We let $u\in k[x]$ be
the monic polynomial with simple zeros in the supersingular points,
and no other zeros. Put
\begin{equation}\label{Qeq}
Q=\prod_{i=1}^{r-1}(x-x_i)^{1+a_i-\nu_i},
\end{equation}
where $a_i$ and $\nu_i$ are defined by (\ref{aidefeq}), i.e.\
\[
\sigma_i=\frac{a_i}{p-1}+\nu_i.
\]

The definition of $\sigma_i$ implies that there exists an $\epsilon\in
k^\times$ such that
\begin{equation}\label{omegaeq}
\omega=\epsilon z\prod_{i\neq r}
(x-x_i)^{\nu_i-1}\,{\rm d}x=\epsilon
\frac{z^p}{Qu^2}\,{\rm d}x.
\end{equation}
Here $z$ is as in (\ref{Kummereq}). 

\begin{lem}\label{corrlem}
 Let $\omega$ be given by (\ref{omegaeq}), and let $u\in k[x]$ be
the monic polynomial with simple zeros exactly in the supersingular points.
\begin{itemize}
\item[(a)] The differential form $\omega$ is 
logarithmic if and only if
\[
D^{p-1}\frac{1}{Qu^2}=-\epsilon^{p-1}\prod_{i=1}
^r(x-x_i)^{p(\nu_i-1)}.
\]
\item[(b)] 
 The differential form $\omega$ is 
exact if and only if
\[
D^{p-1}\frac{1}{Qu^2}=0.
\]
\item[(c)] Suppose that 
$(Z, \omega)$ is a deformation datum. Then $\Res_x \frac{1}{Qu^2}=0,$ 
{for $x$ supersingular}.
\end{itemize}
\end{lem}

\proof It is well known
that $\omega=F\, {\rm d}x$ is logarithmic if and only if
$D^{p-1}F=-F^p$. (For an outline of the proof see \cite[Exercise
9.6]{GS}). Therefore (\ref{Qeq}) implies that $\omega$ is logarithmic
if and only if
\begin{equation}\label{Cartiereq}
D^{p-1}\frac{1}{Qu^2}=-\epsilon^{p-1}\prod_{i\neq r}(x-x_i)^{p(\nu_i-1)}.
\end{equation}
This implies (a).  It is well known that $\omega$ is exact if and only
if $D^{p-1}F=z D^{p-1} (1/Qu^2)=0$. This proves (b).

Let $\tau\in X$ be a supersingular point  and write
\[
   \frac{1}{Qu^2} \;=\; \sum_{n\geq -2} c_n (x-\tau)^n.
\]
Then (\ref{Cartiereq}) implies that
\[  
   D^{p-1}\frac{1}{Qu^2} = -[\frac{c_{-1}}{(x-\tau)^p}+
        c_{p-1}+\cdots] =
        -\epsilon^{p-1}\prod_{i\neq r}(x-x_i)^{p(\nu_i-1)}.
\]
Here $\epsilon\neq 0$ if and only if $\omega$ is logarithmic.  We
conclude that $c_{-1}=0$, since $\tau\neq x_i$ for some $i\in\{1,
\ldots, r\}$.  This proves (b).  \Endproof

 The following proof is inspired by
\cite[Lemma 3]{Beukers}. A similar argument can be found in
\cite[Theorem 5]{Ihara74}. A special case of the result can be found
in \cite[Proposition 3.2]{mcav}.

\begin{prop}\label{corrprop1}
Let $(Z, \omega)$ be a deformation datum and $u$ be the monic
polynomial with simple zeros exactly in the supersingular points.
\begin{itemize}
\item[(a)]  There exists a polynomial
$P_2=d(d+a_r)x^{r-2}+\beta_{r-3}x^{r-3}+\cdots +\beta_0\in k[x]$ such that
\[
L_\omega(u):=u''+\frac{Q'}{Q}u'+\frac{P_2}{\prod_{i=1}^{r-1}(x-x_i)}u=0.
\]
\item[(b)] The differential operator $L_\omega$ defined in (a) is
normalized. It has singularities at $x_1, \ldots, x_r=\infty$. Its
local exponents at $x_i\neq \infty$ (resp.\ $x_r=\infty$) are
$\{-\sigma_i ; 0\}\pmod{p}$ (resp.\ $\{-\deg(u); -\deg(u)
+\sigma_r\}
\pmod{p}$).
\item[(c)] The $p$-curvature of $L$ is nilpotent but nonzero.
\end{itemize}
\end{prop}

\proof
Let $\tau\in X$ be a supersingular point. Lemma \ref{corrlem}.(c) implies that
\[
\Res_{x=\tau}\frac{1}{Qu^2}=0.
\]
The proof of \cite[Proposition 3.2]{mcav} also applies in our more
general situation. We deduce that $u$ is a solution to a Fuchsian
differential equation as in (a).

Since $p_2:=P_2/\prod_{i=1}^{r-1} (x-x_i)$ has at most simple poles,
it follows that $L_\omega$ is normalized. The singularities and local
exponents are easily read off from the explicit expression for
$L_\omega$. This proves (b).  Part (c) follows from \cite[Proposition
4.8]{indi}. \Endproof

The following proposition gives a converse to Proposition
\ref{corrprop1}.  It is a simplified version  of \cite[Proposition 
5.3]{indi} which is stated in a different language.

\begin{prop}\label{existenceprop}
  Let $L=(\partial/\partial t)^2+p_1(\partial/\partial t)+p_2$ be a
  normalized second order differential operator with regular
  singularities in $x_1, \ldots, x_r=\infty$ and local exponents
  $\{\alpha_i, 0\}$ (resp.\ $\{-d , -d+\alpha_r\}$) at $x_i$ for $i\neq r$
  (resp.\ $x_r$). Suppose that $L$ has nilpotent but nonzero
  $p$-curvature. Let $u$ be a polynomial solution of minimal
  degree.   Then the pair $(L, u)$ is associated to a 
  deformation datum $(Z, \omega)$ via the construction of Proposition
  \ref{corrprop1}.
\end{prop}

\proof Let $L$ be as in the statement of the proposition, i.e.\ $p_1$
and $p_2$ are as in (\ref{denormaleq}). We define 
\[
Q=\prod_{i=1}^{r-1} (x-x_i)^{[1-\alpha_i]}.
\]
Here $[a]$ denotes the unique integer satisfying $0\leq [a]<p$ and
$[a]\equiv a\pmod{p}$.  We have $p_1=Q'/Q$.

It follows from Proposition \ref{pcurvprop} that
\[
D^{p-1}\frac{1}{Qu^2}\neq 0.
\]
 Denote by
$x_{r+1}, \ldots x_s$ the points of $X$, different from $x_1, \ldots
x_r$, such that
\[
\ord_{x_i} D^{p-1}
\frac{1}{Qu^2}\neq 0.
\]
For $i=1, \ldots, s$  with $i\neq r$ we define nonnegative integers $\nu_i$ by
\[
\nu_i=
\frac{1}{p}\left(\ord_{x_i}D^{p-1}\frac{1}{Qu^2}\right)+1.
\] 
Moreover, we define integers $0\leq a_i<p$ by $a_i\equiv
-\alpha_i+\nu_i \pmod{p}$, and put $\sigma_i =a_i/(p-1)+\nu_i$.  For
$i=r$ we define $\sigma_r$ by the relation of Lemma
\ref{defodatlem}.(d), and $a_r, \nu_r$ by (\ref{aidefeq}). Note that
$\sigma_i\equiv \alpha_i\pmod{p}$ if $i\leq r$ and $\sigma_i\equiv
1\pmod{p}$ for $i>s$.

After replacing $k$ by a larger algebraically closed field, we may
suppose that all points $x_i$ are rational over $k$.  We let $Z/k$ be the
smooth projective curve defined by the Kummer equation
(\ref{Kummereq}). Lemma \ref{corrlem}.(a) implies that there exists an
$\epsilon \in \kappa^\times$ such that the differential form $\omega$
defined by (\ref{omegaeq}) is logarithmic.

We claim that $\sigma_r\geq 0$. Indeed, if $\sigma_r<0$ then the
differential form $\omega$ has a pole of order strictly larger than
$1$ at $x=x_r=\infty$. But this is impossible since $\omega$ is a
logarithmic differential form.  This implies that $(Z, \omega)$
defines a deformation datum. Its signature is $(\sigma_i)$.
\Endproof

\begin{rem}\label{Dworkrem}
Proposition \ref{Dworkdegprop} follows from Proposition
\ref{existenceprop}. Namely, let $L$ be a differential operator with
nilpotent but nonzero $p$-curvature and let $u$ be a polynomial
solution of minimal degree. To prove Proposition \ref{Dworkdegprop},
we may assume that $L$ is normalized (\S \ref{normalsec}). Proposition
\ref{existenceprop} implies that $(L, u)$ corresponds to a deformation
datum $(Z, \omega)$. Therefore Lemma \ref{defodatlem}.(d) together
with the estimates $\sigma_i\geq t_i/(p-1)$ for $i\leq r$ and
$\sigma_i>1$ for $i>r$ implies that
\[
2\deg(u)\leq (r-2)(p-1)-\sum_{i=1}^r t_i.
\]
\end{rem}

\section{Deformations of $\mup$-torsors}\label{mupsec}
In \S \ref{accsec} we define the moduli space $\N_{0,r}({\boldsymbol
\alpha})$ of differential operators $L$ with nilpotent but nonzero
$p$-curvature and local exponents ${\boldsymbol \alpha}$. The results
of \S \ref{corrsec} imply that $\N_{0,r}({\boldsymbol \alpha})$ also
parameterizes deformation data. We are interested in the dimension of
the irreducible components of $\N_{0,r}({\boldsymbol \alpha})$.  Our
main tool is a result of Wewers (\cite{cotang}) on the deformation of
deformation data. Below we see that this may be translated into the
deformation of $\mup$-torsors.  In this section, we adapt some results
of Wewers to our situation.

Let $(Z, \omega)$ a deformation datum of type $(H, \chi)$.
 We denote by $\G$ the group scheme $\mup\rtimes_\chi H$, as defined
 in \cite[Section 4.1]{cotang}. We associate to the deformation datum
 $(Z, \omega)$ a singular curve $Y$ together with an action of the
 group scheme $\G$, as in
\cite[Construction 4.3]{cotang}. Since $\omega$ is a logarithmic
differential form, locally on $Z$ it may be written as 
\[
\omega=\frac{{\rm d} f}{f},
\]
for some meromorphic functions $f$ on $Z$.
Define $Y$, locally on $Z$, by the equation $y^p=f$. 
Then $\G$ obviously acts on $Y$ and the natural map
$Y\to X$ is a $\G$-torsor outside the branch points of
$Z\to X$ (\cite[Remark 4.6.i]{cotang}). Moreover, \cite[Remark
  4.6.ii]{cotang} implies that $Y$ is
generically smooth.

Let $\Cf_k$ be the category of local artinian $k$-algebras of equal
characteristic $p$. 
A {\sl $\G$-equivariant deformation} of $Y$ over an object $A$
of $\Cf_k$ is a flat $R$-scheme $Y_R$ together with an action of $\G$
and an $\G$-equivariant isomorphism $Y\simeq Y_R\otimes_R
k$. We consider  the deformation functor
\[
R\mapsto \Def(Y, \G)(R)
\]
which sends $R\in \Cf_k$ to the set of isomorphism classes of
$\G$-equivariant deformations of $Y$ over $R$.
Let 
\[
R\mapsto \Def(X; x_i\, |\, i\in \BB)(R)
\]
be the deformation functor which sends $R$ to the set of isomorphism
classes of deformations of the pointed curve $(X; x_i\,|\, i\in
\BB)$. We consider the points $x_i$ on $X$ to be ordered. Moreover, we
consider the $x_i$ up to the action of $\PGL_2(k)$, i.e.\ we suppose
that $x_1=0, x_2=1, x_r=\infty$.  We obtain a natural transformation
\begin{equation}\label{defoeq}
\Def(Y, \G)\To \Def(X; x_i\, |\, i\in \BB).
\end{equation}

\begin{prop}\label{defosmoothprop}
  The deformation functor $\Def(Y, \G)$ is formally smooth. 
\end{prop}

\proof We use the terminology of \cite{cotang}. The proposition follows from
\cite[Theorem 4.8]{cotang}, if we show that $\EE {\rm
  xt}^2_{\G}(\LL_{Y/k}, \OO_{Y})=0$.

In our situation, the integer $s=\dim_{\FF_p}V$ of \cite{cotang} equals one.
 This
implies that the sheaf $\E xt_{\G}^1(\LL_{Y/k}, \OO_{Y})$
has support in isolated points (namely the critical points of the deformation
datum). Since $H^1(X, \E xt_{\G}^1(\LL_{Y/k},
\OO_{Y}))=0$, it follows from \cite[(43)]{cotang} that
\[
\EE{\rm xt}_{\G}^2(\LL_{Y/k}, \OO_Y)=0.
\]
This implies that the deformation problem is formally smooth.
\Endproof

For every $i\in \BB$, we let $\hat{Y}_i$ be the completion of $Y$ at
$x_i$.  Let $R\in \Cf_k$ and $Y_R$ be a $\G$-equivariant deformation
of $Y$. Write $(g_{R}:Z_R\to X_R, \omega_R)$ for the corresponding
deformation datum. Let $i\in \BB$  and
choose a point $z_i$ of $\Zb_0$ above $x_i\in X$. Let $H_i\subset H$
be the decomposition group of $z_i$. There exists a local parameter
$t=t_i$ of $z_i$ on $Z_R$ and a character $\chi_i:H_i\to R^\times$
such that $\OO_{Z_R, z_i}=R[[t]]$ and $h^\ast t_i=\chi_i(h)\cdot
t_i$ for all $h\in H_i$. We denote by $\hat{Y}_{i,R}$ the completion
of $Y_R$ at $x_i$; this is an equivariant deformation of
$\hat{Y}_i$. We obtain a morphism
\begin{equation}\label{locgleq}
{\rm locgl}:\Def(Y,  \G)\To \prod_{i\in \BB} \Def(\hat{Y}_i,
\G)
\end{equation}
called the {\sl local-global morphism} (\cite[Section 5.3]{cotang}).
The $\G$-equivariant deformation $Y_R$ is
called {\sl locally trivial} if it lies in the kernel of the local global morphism.
We denote by
\[
\Def(Y, \G)^{\lt}\subset \Def(Y, \G)
\]
the subfunctor parameterizing locally trivial deformations; this is
the image of the local-global morphism. We write
\[
\Def(Y, \G)^{\lt}=\prod_{i\in\BB} \Def(\hat{Y}_i, \G)^\dagger.
\]

\begin{lem}\label{ltdimlem}
  The tangent space to  $\Def(Y, \G)^{\lt}$ has  dimension 
\[
\frac{1}{p-1}(\sum_{i\in \BB} a_i)-1=N.
\]
\end{lem}

\proof   The tangent
space to the deformation functor $\Def(Y, \G)^{\lt}$ is
\[
H^1(X, \Hom(\LL_{Y/k}, \OO_{Y}))=H^1(X,
\M^{H}),
\]
\cite[Proposition 4.10]{cotang}. Here $\M^{H}$ is defined in \cite[Section
4.3]{cotang}. In our situation it is the sheaf of derivations $D$ of
$\OO_{X}$ such that $D(\omega)$ is a regular function on $Z$. The
proof of \cite[Lemma 5.3]{cotang} implies that $\M^{H}$ is isomorphic to
$((g)_\ast \OO_{Z})_\chi$. A local calculation shows that
\[
\deg(\M^{H})=-\sum_{i\in \BB}\frac{a_i}{p-1}=-N-1.
\]
  By the Riemann--Roch Theorem, the dimension of
$H^1(X, \M^{H})$ equals $-1+(\sum_{i\in \BB} a_i)/(p-1)$. This
proves the lemma.  \Endproof

To compute the dimension of certain components of the moduli space of
deformation data (\S \ref{accsec}), we need to modify the concept of
locally trivial deformations, as defined in \cite{cotang}. The reason
for this is that we allow critical points with ramification invariant
$\sigma_x>2$. This happens, for example,  for the spikes.
 In \cite{cotang}, Wewers focused on {\sl special deformation
data}, which satisfy $\sigma_x\leq 2$. 

Let $x_i$ be a critical point of the deformation datum $(Z, \omega)$,
and let $\sigma_i=a_i/(p-1)+\nu_i$ be its ramification invariant.  In
our situation, the local-global morphism (\ref{locgleq}) is formally
smooth (\cite[Remark 5.12]{cotang}). 
Therefore the following
proposition is proved like \cite[Theorem 5.11.(i)]{cotang}.

\begin{prop}\label{locglprop}
The functor $\Def(\hat{Y}_i, \G)^\dagger$ admits a versal deformation
over the ring
\[
\tilde{R}_i=W(k)[[t_{i, 0}, \ldots, t_{i, \nu_i-1}]].
\]
\end{prop}

\begin{defn}\label{Rastdef}
Let $x_i$ be a critical point of the deformation datum and write
$I_i=\{0\leq j\leq \nu_i\,|\, j\equiv a_i\pmod{p}\}$. Put
$\epsilon_i=|I_i|$. 
We define
\[
\tilde{R}_i^\ast=\tilde{R}_i/\langle t_{i, j}, j\not\in I_i\rangle.
\]
\end{defn}

We now give an interpretation of the rings $\tilde{R}_i^\ast$. Let
$v_i$ be a local parameter of a point $z_i$ of $Z$ above $x_i$. Put
$n_i=(p-1)/\gcd(p-1, a_i)$. The definition of the integer $a_i$
implies that
\[
(\varphi^\ast)^{n_i} v_i\equiv v_i^{b_i} \pmod{v_i^2}, 
\]
where $b_i a_i/n_i\equiv 1\pmod{(p-1)/n_i}$.  Locally around $x_i$, we
 may choose a function $f$ such that $\hat{Y}_i$ is given by the
 Kummer equation $y^p=f$, where $f=1+v_i^{(a_i+(p-1)\nu_i)/n_i}$,
 since $\omega:={\rm d}f/f$ has a zero of order
 $(a_i+(p-1)\nu_i)/n_i-1$ at $z_i$.

 Let
$\hat{Y}^\ast_i$ denote the restriction of the universal deformation
of $\hat{Y}_i$ to $\Spec(\tilde{R}_i^\ast)$. Locally around $x_i$, the germ
$\hat{Y}_i$ is given by a Kummer equation
\begin{equation}\label{Kummereq2}
y^p=F, \qquad \mbox{where }F=f+\sum_{j=0}^{\nu_i-1} t_{i,
j}z^{(a_i+(p-1)j)/n_j},
\end{equation}
since $H$ acts on $\hat{Y}_i$. The differential form corresponding to
$\hat{Y}_i$ is given by ${\rm d}F/F$.

 We claim that $\Spec(\tilde{R}^\ast_i)\subset \Spec(\tilde{R}_i)$ is
exactly the locus such that the restriction of ${\rm d}F/F$ has a zero
of order $(a_i+(p-1)\nu_i)/n_i-1$ at $z_i$, i.e.\ such that the
ramification invariant remains constant. Namely $\partial F/\partial
v_i$ has a single zero in $v_i$ if and only if $t_{i, j}=0$ for all
$j$ such that $a_i-j\not\equiv 0\pmod{p}$. This proves the claim.

Let $\Def(Y, \G; {\boldsymbol \sigma})\subset \Def(Y, \G)$ be the
deformation problem of $\G$-equivariant deformations with fixed
signature. For every critical point of the deformation datum, we
denote by $\Def(\hat{Y}_i, \G; \sigma_i)\subset \Def(\hat{Y}_i, \G)$
the subfunctor parameterizing local deformations with given signature.

\begin{thm}\label{sigmadimthm}
\begin{itemize}
\item[(a)] The deformation problem $\Def(Y, \G; {\boldsymbol \sigma})$
is formally smooth.
\item[(b)] The functor $\Def(\hat{Y}_i, \G; \sigma_i)$ admits a versal
deformation over the ring $\tilde{R}_i^\ast$.
\item[(c)]  The dimension of $\Def(Y, \G; {\boldsymbol \sigma})$ is 
\[
N+\sum_{i\in \BB} \epsilon_i.
\]
\end{itemize}
\end{thm}

\proof Part (a) follows from Proposition \ref{defosmoothprop}. Part (b)
follows from the above discussion.

We have already remarked that the local-global morphism
(\ref{locgleq}) is formally smooth. This implies that the dimension of
the deformation problem $\Def(Y, \G; {\boldsymbol \sigma})$ is equal to the
dimension of its tangent space. It follows that
\[
\begin{split}
\dim \Def(Y, \G; {\boldsymbol \sigma})&=\dim \Def(Y, \G)^\lt+\sum_{i\in
\BB}(\dim \Def(\hat{Y}_i, \G)-\dim \Def(\hat{Y}_i, \G,
\sigma_i)\\
&=N+\sum_{i\in \BB} \epsilon_i.
\end{split}
\]
The last equality follows from Lemma \ref{ltdimlem} together with 
(a).  This proves the theorem.  \Endproof

The following lemma gives a bound on $\epsilon_i$. We assume that
$p>r-2$; this will be assumed in \S \ref{accsec},
as well.

\begin{lem}\label{epsilonboundlem}
Assume that $p>r-2$.
\begin{itemize}
\item[(a)] We have that $\epsilon_i\in \{0, 1\}$.
\item[(b)] If $x=x_i$ is supersingular, then $\epsilon_i=0$.
\item[(c)] If $x=x_i$ is a spike, then $\epsilon_i=1$.
\end{itemize}
\end{lem}

\proof
Lemma \ref{defodatlem}.(d) implies that
\begin{equation}\label{vceq}
\sum_{i=1}^r \sigma_i+\frac{2d}{p-1}+\sum_{i=r+1}^s (\sigma_i-1)=r-2.
\end{equation}
Here $d$ is the number of supersingular points.
This implies that if $x=x_i$ is a singularity (i.e. $1\leq i\leq r$),
then $\sigma_i\leq r-2<p$. We conclude that $\epsilon_i\leq 1$. Since
the supersingular points have ramification invariant
$\sigma=(p+1)/(p-1)$, it follows that they have $\epsilon=0$.

Now let $x=x_i$ be a spike. Then $\sigma_i\leq r-2+1=r-1\leq
p$. Moreover, $\sigma_i\equiv 1\pmod{p}$. This implies that
$\nu_i\equiv 1+a_i\bmod{p}$. Since $\sigma_i\neq 1$, it follows that
$2\leq \nu_i<p$. This implies that $I_i=\{ \nu_i-1\}$ and hence that
$\epsilon_i=1$.
\Endproof

The following proposition characterizes the signatures ${\boldsymbol
\sigma}$ for which the dimension of $\Def(Y, \G; {\boldsymbol
\sigma})$ is maximal.

\begin{prop}\label{dimmaxprop}
\begin{itemize}
\item[(a)] Let $\pi_{\boldsymbol \sigma}:\Def(Y, \G; {\boldsymbol
\sigma})\to \Def(X; \{x_1, \ldots, x_r\})$ be the natural map which
sends a deformation datum to its set of singularities. Then
$\pi_{\boldsymbol \sigma}$ is finite.
\item[(b)] Suppose that $p>r-2$. Then the dimension of $\Def(Y, \G;
{\boldsymbol \sigma})$ equals $r-3=\dim \Def(X; \{x_1, \ldots, x_r\})$
if and only if
\[
\nu_i=\begin{cases} 0&\mbox{ for } 1\leq i\leq r,\\
2&\mbox{ for } r+1\leq i\leq s.
\end{cases}
\]
\end{itemize}
\end{prop}

\proof Part (a) follows from the main result of \cite{Dwork}, by using
the correspondence between deformation data and differential operators
(\S \ref{corrsec}). (Compare to Proposition \ref{Mochizukiprop}.(a).)

Suppose now that $p>r-2$.
We deduce from Theorem \ref{sigmadimthm} and (\ref{vceq}) that the
dimension of $\Def(Y, \G; {\boldsymbol \sigma})$ equals $r-3$ if and only if
\[
\sum_{i=1}^r (\nu_i-\epsilon_i)+\sum_{i=r+1}^s(\nu_i-\epsilon_i-1)=0.
\]
Let $x=x_i$ be a singularity, i.e.\ $1\leq i\leq r$, and suppose that
$\epsilon_i=1$. Definition \ref{ddssdef} implies that
$\sigma_i\not\equiv 1\pmod{p}$. It follows therefore that $\nu_i\geq
3$, and hence that $\nu_i-\epsilon_i>0$. If $x=x_i$ is a spike, i.e.\
$i>r$, we deduce from Lemma \ref{epsilonboundlem}.(c) that
$\nu_i-\epsilon_i-1=\nu_i-2$. This proves the proposition.
\Endproof

\section{The accessory parameter problem}\label{accsec}
This section contains the main results of the paper. In \S
\ref{accintrosec} we introduce the stacks $\N_{0,r}({\boldsymbol
\alpha})$ parameterizing differential equations with nilpotent but
nonzero $p$-curvature. We also define the strength and prove a
necessary condition for the substack $\N_{0,r}({\boldsymbol
\alpha};n)$ of operators of strength $n$ to be nonempty (Lemma
\ref{nlem}). In \S \ref{dimsec}, we prove a result on the dimension of
irreducible components of $\N_{0,r}({\boldsymbol \alpha};n)$ for those
$({\boldsymbol \alpha}, n)$ for which this stack is nonempty, using
the results of \S \ref{mupsec}. In \S \ref{logmonsec} we consider the
case of logarithmic local monodromy.

\subsection{Statement of the problem}\label{accintrosec}
Let $r\geq 3$ be an integer.  We denote by $\M_{0,r}/\FF_p$ the stack
parameterizing $r$-marked curves $(X; \{x_1, \ldots, x_r\})$ of genus
zero.  Fix a set ${\boldsymbol \alpha}=(\alpha_1, \ldots, \alpha)$
with $\alpha_i\in \FF_p$.

Dwork's {\sl accessory parameter problem} asks to determine the space
of all normalized differential operators $L$ with nilpotent
$p$-curvature. This amounts the find all $(x_i;\beta_i)$ such that the
differential operator given by (\ref{denormaleq}) is normalized and
admits a polynomial solution. Dwork (\cite{Dwork}) proves that the
algebraic space, $V_N$, of such differential operators is a complete
intersection. Moreover, Dwork shows that the natural projection of
$V_N$ on  $\M_{0,r}$  has degree $p^{r-3}$
(\cite[Corollary 4.2]{Dwork}).

In this section, we are interested in the moduli space of normalized
differential operators with nilpotent but nonzero $p$-curvature, or,
equivalently, the moduli space of deformation data.

\begin{prop}\label{Mochizukiprop}
\begin{itemize}
\item[(a)] 
 There exists a stack $\N_{0,r}({\boldsymbol \alpha})$
parameterizing normalized differential operators with nilpotent,
nonzero $p$-curvature with local exponents ${\boldsymbol \alpha}$.
\item[(b)] Let $\pi:\N_{0,r}({\boldsymbol \alpha})\to \M_{0,r}$ be the
natural projection which sends $L$ to its set of singularities. Then
$\pi$ is finite. 
\item[(c)] The degree of $\pi$ is less than or equal to $p^{r-3}$,
with equality if and only if there do not exist normalized
differential operators with local exponents ${\boldsymbol \alpha}$
with zero $p$-curvature.
\end{itemize}
\end{prop}

\proof Mochizuki (\cite[Chapter IV]{Mochizuki2}) proves the existence
of a stack parameterizing indigenous bundles with local exponents
${\boldsymbol \alpha}$.  It follows from \cite[Theorem 4.11]{indi} and
Proposition \ref{existenceprop} that there is equivalence between
indigenous bundles and normalized differential operators with
nilpotent, nonzero $p$-curvature.
Parts (b) and (c)  follow, for example, from \cite[Corollary 4.2]{Dwork}.
 \Endproof

  Let $L$ be a
 normalized differential operator with nilpotent but nonzero
 $p$-curvature. Then Proposition \ref{existenceprop} implies that $L$
 defines a deformation datum $(Z, \omega)$. Therefore we may define
 the signature ${\boldsymbol \sigma}=(\sigma_i)_{i\in \BB}$, as in the
 proof of Proposition \ref{existenceprop}.  In particular, we may
 define the spikes $x_{r+1}, \ldots, x_s$. 

\begin{defn}\label{nxdef}
Let $x$ be a critical point of the deformation datum corresponding
to $L$. Define
\[
n_x=\begin{cases} 0&\mbox{ if } x \mbox{ is supersingular},\\
(p-1)\sigma_x&\mbox{ if } x \mbox{ is a singularity},\\
(p-1)(\sigma_x-1)&\mbox{ if } x \mbox{ is a spike}.
\end{cases}
\]
The number 
$n:=\sum_x n_x$ is called the {\sl strength} of $L$.
\end{defn}

Definition \ref{ddssdef} implies that $n_x$ is a nonnegative integer.
We refer to \cite{indi} for an explanation of this notion; suitably
defined, $n_x$ is the order of vanishing of the $p$-curvature at $x$.
The proof of Proposition \ref{existenceprop} implies that
\[
n_x=\ord_x\prod_{i=1}^{r-1}(x-x_i)^p D^{p-1}\frac{1}{Qu^2}.
\]
This gives a concrete interpretation of $n_x$.
The terminology `strength' was introduced by Mochizuki
(\cite[Introduction, \S 1.2]{Mochizuki2}). Note that the index $d$ in
\cite{Mochizuki2} denotes the strength and not the number of
supersingular points. The following lemma states a few properties of
$n$. The lemma follows immediately from (\ref{vceq}).

\begin{lem}\label{nlem}
Suppose that $n$ is the strength of a deformation datum. Then
\begin{itemize}
\item[(a)] $n= (r-2)(p-1)-2d$,
\item[(b)] $n$ is
an even integer which is equivalent to $-\sum_{i=1}^r
\alpha_i\pmod{p}$,
\item[(c)] $n\geq \sum_{\alpha_i\neq 0}(p-\alpha_i)$.
\end{itemize}
\end{lem}

We denote by $\N_{0,r}({\boldsymbol \alpha}; n)$ the substack of
 $\N_{0,r}({\boldsymbol \alpha})$ parameterizing normalized
 differential operators with nilpotent and nonzero $p$-curvature and
 strength $n$.  The {\sl strong accessory parameter problem} asks for
 the structure of $\N_{0,r}({\boldsymbol \alpha}; n)$. For example, it
 is natural to ask for which $n$ the stack $\N_{0,r}({\boldsymbol
 \alpha}; n)$ is nonempty or has maximal dimension (i.e.\
 $\dim\N_{0,r}({\boldsymbol \alpha}; n)=\dim \M_{0,r}=r-3$). Also, one
 would like to know the degree of $\N_{0,r}({\boldsymbol \alpha};
 n)\to \M_{0,r}$, for those $n$ for which $\N_{0,r}[n]$ has maximal
 dimension.

\begin{rem} Lemma \ref{nlem} gives a necessary condition for the nonemptyness
 of $\N_{0,r}({\boldsymbol \alpha}; n)$. Namely for given
${\boldsymbol \alpha}$, a necessary condition for
$\N_{0,r}({\boldsymbol \alpha})$ to be nonempty is that there exists a
nonnegative integer $n$ such that the (in)equalities of Lemma
\ref{nlem} are satisfied. We leave it to the reader to check that this
is indeed a nontrivial condition.
\end{rem}

As far as I know, the only results in this direction are for $r\leq
4$. These results are summarized in Example \ref{Gaussexa}.

\begin{exa}\label{Gaussexa} 
Suppose that $\alpha_i=0$ for $i=1, \ldots, r$. This case we consider
in more detail in \S \ref{logmonsec}.

  (a) If $r=3$ it follows from Proposition \ref{pcurv0prop} that $d<
p$. Since $d\equiv -1/2\pmod{p}$, it follows that $d=(p-1)/2$ and hence that
$\N_{0,3}({\boldsymbol \alpha})=\N_{0,3}({\boldsymbol \alpha};0)$. In
fact, this space consists of one point. Let $r=3$ and suppose that
$x_1=0, x_2=1, x_3=\infty$. Then the corresponding differential
operator $L$ is the Gauss hypergeometric differential operator
(\cite[Example 4.5]{indi}):
\[
L=(\partial/\partial t)^2+\frac{2t-1}{t(t-1)}(\partial/\partial
t)+\frac{1}{4t(t-1)}.
\]
Its monic polynomial solution of minimal degree is the Hasse invariant:
\[
\Phi=\sum_{i=0}^{(p-1)/2}\binom{p-1}{i}^2 t^i.
\]
The corresponding deformation datum describes the stable reduction of
the cover of modular curves $X(2p)\to X(2)$ (\cite{mcav}).

(b) Similarly, if $r=4$ it follows from Proposition \ref{Dworkdegprop}
or from \cite{Beukers} that $d=p-1$, hence $\N_{0,4}({\boldsymbol
\alpha})=\N_{0,4}({\boldsymbol \alpha}; 0)$. Lemma
\ref{pcurvnonzerolem} and \cite[Corollary 4.2]{Dwork} imply, moreover,
that $\N_{0,4}({\boldsymbol \alpha})$ is nonempty.
 In \cite[\S 6.2]{indi} it is shown that
there exists a connected component which has degree one over
$\M_{0,4}$.
\end{exa}

\subsection{The dimension of the components of $\N_{0,r}({\boldsymbol
\alpha};n)$}
\label{dimsec}
 Suppose that $\N:=\N_{0,r}({\boldsymbol \alpha}; n)$ is
nonempty.  In this section, we characterize the those irreducible
components of $\N$ which have maximal dimension $r-3$.  This is a
direct consequence of the results of \S \ref{mupsec}.

We fix local exponents ${\boldsymbol \alpha}=(\alpha_1, \ldots,
\alpha_r)$ and a nonnegative integer $n\equiv
-\sum_{i=1}^r\alpha_i\pmod{p}$ such that $d:=[(r-2)(p-1)-n]/2$ is a
nonnegative integer and $n\geq \sum_{\alpha_i\neq 0} (p-\alpha_i)$
(Lemma \ref{nlem}). Let $s$ be defined by 
\[
s-r=\frac{1}{p}[(r-2)(p-1)-2d-\sum_{i:\, \alpha_i\neq 0} (p-\alpha_i)].
\]
It follows from the assumptions on $n$ that $s$ is a nonnegative
integer.

Let $\BB=\{1, \ldots, d+s\}$.
We now define a signature ${\boldsymbol \sigma}={\boldsymbol
\sigma_{\boldsymbol \alpha}}=(\sigma_i)_{i\in \BB}$ as follows.  For
$1\leq i\leq r$, we put
\[
\sigma_i=\begin{cases} 0&\mbox{ if }\alpha_i=0,\\
(p-\alpha_i)/(p-1)&\mbox{ otherwise}.
\end{cases}
\] 
For $r+1\leq i\leq s$, we put $\sigma_i=(2p-1)/(p-1)$. For $i>s$, we
put $\sigma_i=(p+1)/(p-1)$. Note that ${\boldsymbol \sigma}$ has $s-r$
spikes and $d$ supersingular points.

\begin{prop}\label{dimprop} Suppose that $p>r-2$ and  that 
$\N_{0,r}({\boldsymbol \alpha}; n)$ is nonempty. Let $\N$ be an
irreducible component of $\N_{0,r}({\boldsymbol \alpha}; n)$ of
dimension $r-3$.  Then there exists a point $L\in \N$ whose signature
is ${\boldsymbol \sigma}_{\boldsymbol \alpha}$.
\end{prop}

\proof Let $\N$ be in the statement of the proposition, and suppose
that $\dim \N=r-3$. Let $L$ be the differential operator corresponding
to the generic point of $\N$. Proposition \ref{dimmaxprop}.(b) implies
that the signature of $L$ satisfies $\sigma_i<1$ for $1\leq i\leq r$
and $\sigma_i=(2p-1)/(p-1)$ if $x_i$ is a spike. This implies that the
signature of $L$ is ${\boldsymbol \sigma}_{\boldsymbol \alpha}$.
\Endproof

\subsection{Examples}\label{exasec}
In this section, we illustrate the results of the previous section
with some examples.

Let $p=7$ and $r=4$. We choose $\alpha_1=\alpha_2=\alpha_3=\alpha_4=5$
and $x_1=0, x_2=1, x_3=\lambda, x_4=\infty$, where we assume that
$\lambda$ is transcendental over $\bar{\FF}_p$. In other words, we
want to determine components of $\N_{0,4}({\boldsymbol \alpha})$ with
dimension $r-3=1$. Put $k=\bar{\FF}_p(\!(\lambda)\!)$.  We consider
differential operators $L$ with singularities in $(x_i)$ and local
exponents ${\boldsymbol \alpha}=(\alpha_i)$, i.e.\ $L$ is given by
(\ref{denormaleq}). However, we do not assume that $L$ is normalized
(\S \ref{normalsec}). Let $\gamma_1, \gamma_2$ be the local exponents
at $\infty$. Therefore it follows from Lemma \ref{deglem} that if
$u\in k[x]$ is a solution of $L$ the $\deg(u)\equiv -\gamma_1\equiv 2,
-\gamma_2\equiv 4\pmod{p}$. (Assuming that $L$ is normalized would
exclude one of these possibilities.)

Set $u=\sum_{i\geq 0} u_ix^i$. One checks that $L(u)=0$ if and only if
the coefficients of $u$ satisfy the recursion
\[
\lambda A_iu_{i+1}=(C_i-\beta)u_i-B_iu_{i-1},
\]
where
\[
\begin{split}
A_i&=(i-1)(i+1-\alpha_1), \qquad B_i=(i+\gamma_1-1)(i+\gamma_2-1),\\
C_i&=i^2(1+\lambda)+i(\lambda(1-\alpha_1-\alpha_2)+1-\alpha_1-\alpha_3).
\end{split}
\]

Propositions \ref{pcurv0prop} and \ref{Dworkdegprop} imply that if $L$
has nilpotent $p$-curvature then $L$ has a solution of degree less
than $p$.

One deduces that $L$ has a solution of degree $2$ if and only if the
accessory parameter $\beta$ satisfies $\beta\in \{0,-1,
-\lambda\}$. Lemma \ref{deglem}.(b) implies in each of these cases
that the solution, $u$, of degree $2$ does not have zeros in
$x=0,1,\lambda$. Moreover, one checks that
\[
\left(\frac{\partial}{\partial
x}\right)^{p-1}\frac{1}{x^3(x-1)^3(x-\lambda)^3u^2}\neq 0,
\]
as rational function of $\lambda$. Therefore it follows from
Proposition \ref{pcurvprop} that the $p$-curvature of $L$ is nonzero
for $\beta\in \{0,-1, -\lambda\}$.  (One could also deduce this from
the recursion for the coefficients of $u$.) Lemmas
\ref{defodatlem}.(d) and \ref{nlem} imply that the signature of $(L,
u)$ is ${\boldsymbol \sigma}=(\sigma_i)$ with $\sigma_i=2/(p-1)=1/3$
for $i=1, \ldots, 4$. In particular, we have no spikes and
$n=2(p-1)-2d=8$.

We conclude that $\N_{0,4}({\boldsymbol \alpha};
8)\subset\N_{0,4}({\boldsymbol \alpha})$ is dense. The degree of
$\pi:\N_{0,4}({\boldsymbol \alpha})\to \M_{0,4}$ is $3$. More
precisely, $ \N_{0,4}({\boldsymbol \alpha}; 8)$ consists of three
irreducible components which each have degree $1$ over
$\M_{0,4}$. 

Similarly, one shows that $L$ has a solution, $u$, of degree $4$ if and only if
\[
5\beta^2+\beta(1+\lambda)+\lambda^2+3\lambda+1=0.
\] 
As before, it follows that $u$ does not have zeros in
$x=0,1,\lambda$. One computes that the residues of
$1/(x^3(x-1)^3(x-\lambda)^3u^2)$ in $x=0,1,\infty$ are zero. This
implies that
\[
\left(\frac{\partial}{\partial
x}\right)^{p-1}\frac{1}{x^3(x-1)^3(x-\lambda)^3u^2}=0
\]
(Lemma \ref{corrlem}.(c)).  It follows that the $p$-curvature of $L$
is zero in this case. In fact, one computes that $L$ also has a
solution of degree $2+p=9$.

\bigskip\noindent
As a second example, we consider $p=13, r=4$ and ${\boldsymbol
\alpha}=(11, 11, 11, 10)$. Let $x_1=0, x_2=1, x_3=\lambda,
x_4=\infty$, with $\lambda$ transcendental over $\bar{\FF}_p$, as in
the previous example. We let $L$ be given by (\ref{denormaleq}), i.e.\
$L$ is a general differential operator with local exponents
${\boldsymbol \alpha}$. Since we do not assume that $L$ is normalized,
it follows that the degree of a polynomial solution of $L$ is
congruent to $1, 4\pmod{p}$. 

Arguing as in the previous example, we
find that $L$ has a solution of degree $1$ if and only if 
\begin{equation}\label{beta1eq}
\beta^2+7(\lambda+1)\beta+\lambda=0.
\end{equation}
Moreover, for $\beta$ satisfying (\ref{beta1eq}), the monic polynomial
solution, $u$, of $L$ does not have a zero in $x=0,1,\lambda$. One
computes that
\[
\left(\frac{\partial}{\partial
x}\right)^{p-1}\frac{1}{x^3(x-1)^3(x-\lambda)^3u^2}
\]
has exactly one zero, which has  order $p$. Moreover, this zero is
not in $x=0,1,\lambda$, since $\lambda$ is transcendental. This also
follows from Proposition \ref{dimprop}. In particular, the
$p$-curvature of $L$ is nonzero. We conclude that $L$ has one spike,
and hence that the signature is $(2/(p-1), 2/(p-1), 2/(p-1), 3/(p-1),
(2p-1)/(p-1))$. It follows that the strength is $n=22$.

We conclude that $\N_{0,4}({\boldsymbol
\alpha};22)\subset\N_{0,4}({\boldsymbol \alpha})$ is dense and that
the degree of $\pi:\N_{0,4}({\boldsymbol \alpha})\to \M_{0,4}$ is $2$.

Similarly, one checks that $L$ has a polynomial solution of degree $4$
if and only if
\begin{equation}\label{beta2eq}
\lambda^3+(2\beta+9)\lambda^2+(9+8\beta^2+4\beta)\lambda+2\beta^3+8\beta^2+
2\beta+1=0.
\end{equation}
One checks that, for $\beta$ satisfying (\ref{beta2eq}), the
$p$-curvature of $L$ is nonvanishing, and that $L$ does not have a
spike. The corresponding signature is therefore ${\boldsymbol
\sigma}=(2/(p-1), 2/(p-1), 2/(p-1), 10/(p-1))$. Hence the strength is
$n=16$. The differential operator is not normalized: following the
convention of \S \ref{normalsec} we have that ${\boldsymbol
\alpha}':=(p-2, p-2, p-2, 3)$ in this case.

We conclude that $\N_{0,4}({\boldsymbol
\alpha}';16)\subset\N_{0,4}({\boldsymbol \alpha}')$ is dense and that
the degree of $\pi:\N_{0,4}({\boldsymbol \alpha}')\to \M_{0,4}$ is $3$.

\subsection{The case of logarithmic local monodromy}\label{logmonsec}
In this section, we consider the strong accessory parameter problem in the
case that all local monodromy matrices are nilpotent.

\begin{defn} Let $L$ be a normalized differential operator with nilpotent 
$p$-curvature, and let ${\boldsymbol \alpha}=(\alpha_i)_{i=1}^r$ be
its local exponents.  We say that $L$ has {\sl logarithmic local
monodromy} if $\alpha_i=0$ for $i=1, \ldots, r$.
\end{defn}

 Let $d$ be a nonnegative integer congruent to $1-r/2 \pmod{p}$ such
that $n:=(p-1)(r-2)-2d$ is nonnegative.  In this section, we drop the
local exponents from the notation, and write $\N_{0,r}[n]$ for the
stack parameterizing differential operators with nilpotent
$p$-curvature and logarithmic local monodromy and strength $n$.  Lemma
\ref{pcurvnonzerolem} implies that the $p$-curvature of a differential
operator is always nonzero, therefore this notation agrees with the
notation in the previous section.

Let $L$ correspond to a point of $ \N_{0,r}[n]$, and let $u$ be a
polynomial solution of $L$ of minimal degree,  $x_1, \ldots,
x_r$  the singularities of $L$. Then $(L, u)$ corresponds to a
deformation datum $(Z, \omega)$ (Proposition
\ref{existenceprop}). Therefore we may use the terminology of \S
\ref{ddgensec}.  Let $x_{r+1}, \ldots, x_s$ be the spikes of $(Z,
\omega)$ and let $\sigma_i=a_i/(p-1)+\nu_i$ be the ramification
invariant of $x_i$ for $1\leq i\leq s$ (compare to the proof of
Proposition \ref{existenceprop}.) It follows from \cite[Proposition
3.6.(i)]{indi} that $\sigma_i=0$ for $i=1, \ldots, r$.

 The following theorem is proved by Mochizuki (\cite[Introduction,
Theorem 1.2]{Mochizuki2}), by using a deformation argument. It is
stronger than Proposition \ref{dimprop}.  Namely, Theorem
\ref{Mdimthm} implies that every differential operator may be deformed
to a differential operator of the same strength such that all spikes
have the minimal possible ramification invariant, namely
$\sigma=(2p-1)/(p-1)$. One expects this to hold for arbitrary local
exponents, as well.

\begin{thm}[Mochizuki]\label{Mdimthm} 
Suppose that $\N_{0,r}[n]$ is nonempty. Then all irreducible components of 
$\N_{0,r}[n]$ have dimension $r-3$.
\end{thm}

The following lemma is a more precise version of Lemma \ref{nlem} in
the case of logarithmic local monodromy.

\begin{lem}\label{slem}
Let $r\geq 3$ and $n\geq 0$ be integers such that $\N_{0,r}[n]$ is
nonempty. 
\begin{itemize}
\item[(a)] Then $n\equiv 0\pmod{2p}$.
\item[(b)]  We have that $ 0\leq n\leq (p-1)(r-2)$.
Moreover,  $r\geq n/p+3$.
\end{itemize}
\end{lem}

\proof The statement that $n\equiv 0\pmod{p}$ follows from Definition
\ref{ddssdef}.(d). Lemma \ref{defodatlem}.(d) implies that
$2d+n=(r-2)(p-1)$, since $\sigma_i=0$ for $i=1,\ldots, r$. It follows
that $n$ is even. 

 The  bounds on $n$  follow from Proposition
\ref{Dworkdegprop}.  The inequality for $r$ follows from (a).\Endproof

To goal of the rest of this section is to prove the following
theorem. This theorem follows immediately from Propositions
\ref{pcuspprop} and  \ref{recursionprop}.

\begin{thm}\label{nonemptythm}
Suppose that $p> r-2$. Let $0\leq n\leq (p-1)(r-2)$ be congruent to $0
\pmod{2p}$.  Then $\N_{0,r}[n]$ is nonempty.
\end{thm}

\begin{prop}\label{pcuspprop}
Suppose that $\N_{0,r}[n]$ is nonempty. Then $\N_{0, r+1}[n]$ is
non\-empty, as well.
\end{prop}

\proof Suppose that $\N_{0,r}[n]$ is nonempty, and let $(L, u)$
correspond to a point of $\N_{0,r}[n]$, i.e.\ $u$ is a polynomial
solution of $L$ of minimal degree $d$. The integer $d$ satisfies
$d=[(r-2)(p-1)-n]/2\leq (p-1)(r-2)/2$. As usual, we denote by $x_1,
\ldots, x_r=\infty$ the singularities of $L$ and by $x_{r+1}, \ldots,
x_s$ the spikes. Proposition \ref{existenceprop} implies that
the pair $(L, u)$ corresponds to a deformation datum $(Z, \omega)$. The
main result of \cite{indi} implies that this deformation datum
corresponds to an indigenous bundle $\E_1$ on  
$(X_1:=\PP^1; x_i)$.

Let $\E_2$ be the indigenous bundle corresponding to the deformation
datum of Example \ref{Gaussexa}.(a); this is the indigenous bundle
corresponding to Gauss' hypergeometric differential equation. The
indigenous bundle $\E_2$ lives on the marked curve $(X_2; \tau_1=0,
\tau_2=1, \tau_3=\infty)$ and has no spikes, i.e.\ $n=0$.

We define a stably marked curve $X$ by identifying the point
$x_r=\infty$ on $X_1$ with the point $\tau_3=\infty$ on
$X_2$. Mochizuki (\cite[\S I.2, page 1008]{Mochizuki1}) shows that the
indigenous bundles $\E_i$ define an indigenous bundle $\E$ on $X$.  We
refer to \cite{Mochizuki1} for the precise definition of an indigenous
bundle on a stably marked curve. The bundle $\E$ is what Mochizuki
calls an {\sl indigenous bundle of restrictable type}. It is shown in
\cite[Proposition 2.11, \S I.2]{Mochizuki1} that $(X, \E)$ deforms to
an indigenous bundle $\tilde{\E}$ on a smooth stably marked curve
$(\tilde{X}; \tilde{x}_i)$ which is a deformation of the stably marked
curve $X$. In particular, the bundle $\tilde{\E}$ (or equivalently,
the corresponding pair $(\tilde{L}, \tilde{u})$) has $r+3-2=r+1$
singular points and is spiked of strength $n$. This proves that
$\N_{0, r+1}[n]$ is nonempty.  \Endproof

One could give an alternative proof of Proposition \ref{pcuspprop}
without using the results of Mochizuki by using the ideas of
\cite[Section 2.5]{Crelle}.

\begin{prop}\label{recursionprop} Suppose that $p> r-2$. Let 
$0\leq n\leq (p-1)(r-2)$ be an integer 
with $n\equiv 0\pmod{2p}$, and let $r=n/p+3$. Then $\N_{0,r}[n]$ is
nonempty.
\end{prop}

\proof
Note that the assumption that $r=n/p+3$ implies that $r$ is odd.
 Example \ref{Gaussexa} implies that the proposition holds for
$r=3$, therefore it is no restriction to suppose that $r\geq 5$, or,
equivalently, $n\geq 2p$. The degree, $d$, of a minimal polynomial
solution is now $d=(p-r+2)/2$.

Let
$\zeta\in k$ be a primitive $(r-2)$th root of unity. It exists since we
assumed that $p>r-2$. We define $x_1=0, x_r=\infty$. Moreover, for
$i=0, \ldots r-3$, we put $x_{i+1}=\zeta^i$. We consider a general
differential operator $L$ with logarithmic local monodromy at $x_1,
\ldots, x_r$. By (\ref{denormaleq}) we have that
\[
L=(x^{r-1}-x)\frac{\partial}{\partial x}^2+
((r-1)x^{r-2}-1)\frac{\partial}{\partial x}+
(d^2x^{r-3}+\beta_{r-4}x^{r-4}+\cdots+\beta_0).
\]
We want to determine $\beta_i$ such that $L$ has a polynomial
solution, $u$, of degree $d$. Write $u=\sum_{i\geq 0} u_ix^i$ with
$u_0=1$. Then setting $L(u)=0$ yields the following recursion for the
coefficients:
\begin{equation}\label{recursioneq}
u_{i-r+3}\left(i-\frac{r-4}{2}\right)^2+\sum_{j=0}^{r-4} u_{i-j} \beta_j=u_{i+1}(i+1)^2.
\end{equation}
This determines the coefficients $u_1, \ldots, u_{p-1}$ uniquely, in
terms of the $\beta_i$, since for $0\leq i\leq p-2$ we have that
$i+1\not\equiv 0\pmod{p}$.

To show the existence of a polynomial solution of degree $d$ of $L$,
we have to show that we may choose $\beta_i\in k$ such that
$u_{d+1}=\cdots u_{d+r-3}=0.$ The recursion (\ref{recursioneq}) then
implies that we may assume that $u_i=0$ for all $i>d$.

One easily deduces from (\ref{recursioneq}) that the total degree of
$u_{i}$ in the variable $\beta_j$ is $i$, and that
\begin{equation}\label{beta0eq}
u_i=\epsilon_i\beta_0^i+\text{  terms of strictly lower degree, for some }
\epsilon_i\in \FF_p^\times.
\end{equation}
 We now homogenize the equations $u_{d+1}=0,
\ldots, u_{d+r-3}=0$, introducing a new variable $\gamma$, and consider
our equations as equations on $\PP^{r-3}$. Proposition \ref{Mochizukiprop}.(b)
implies that the solutions space of these equations is zero
dimensional. Bezout's Theorem implies that the total number of
solutions (counted with multiplicity) is $(d+1)(d+2)\cdots(d+r-3)$. It
remains to show that the number of solutions with $\gamma=0$ is strictly
less than this number.

We now compute the solutions with $\gamma=0$. From the particular form
of the equations, we deduce that if $\gamma=0$ then also
$\beta_0=0$. We now eliminate the variable $\beta_0$ from the
equations $u_{d+2}=0, \ldots, u_{d+r-3}=0$. We first replace
$u_{d+i}=0$ by $\E_{d+i}:=u_{d+i}-\beta_0^{i-1}u_{d+1}\mu_i=0$, where
$\mu_i=\epsilon_{d+i}/\epsilon_{d+1}\in \FF_p^\times$.  Then
$\E_{d+i}$ is divisible by $\gamma$. Dividing out by a power of
$\gamma$ and substituting $\beta_0=0$, we obtain therefore for $i>1$ a
new equation of degree strictly smaller that $d+i$. Since $r-3>1$, we
conclude that not all solutions of our system of equations are on the
hyperplane at $\infty$. This shows that the differential equation has
a solution $u$ of degree less than or equal to $d$. Since $u$ is
nontrivial and $d$ is the smallest possible degree of a solution, we
conclude that $\deg(u)=d$.  This proves the proposition.  \Endproof

\flushright{ Irene I.~Bouw\\
Institut f\"ur reine Mathematik\\
Universit{\"a}t Ulm\\
D-89069 Ulm\\
irene.bouw@uni-ulm.de}

\end{document}